\documentclass[12pt,notitlepage,twoside,a4paper]{amsart}
 \usepackage{amsfonts}

\usepackage{amsmath,amssymb,enumerate}

\usepackage{epsfig,fancyhdr,color}

\usepackage{amssymb}
\usepackage{amsmath,amsthm}
\usepackage{latexsym}
\usepackage{amscd}
\usepackage{psfrag}
\usepackage{graphicx}
\usepackage[latin1]{inputenc}
\usepackage[all]{xy}
\usepackage[mathcal]{eucal}

\definecolor{NoteColor}{rgb}{1,0,0}

\newcommand{\degree}{\ensuremath{^\circ}}

\renewcommand{\textsc}{\textcolor{red}}

\newtheorem*{theorem 1}{\rm\bf Proposition 1}
\newtheorem*{theorem 2}{\rm\bf Proposition 2}

\theoremstyle{definition}

\theoremstyle{remark}

\def\interieur#1{\mathord{\mathop{\kern 0pt #1}\limits^\circ}}

\title[Hyperbolic geometry]{Hyperbolic geometry in the work of J. H. Lambert}

\author{Athanase Papadopoulos}
\address{Athanase Papadopoulos, Institut de Recherche Math\'ematique Avanc\'ee,
Universit\'e de Strasbourg and CNRS,
7 rue Ren\'e Descartes,
 67084 Strasbourg Cedex, France.} 
 \email{athanase.papadopoulos@math.unistra.fr}
\author{Guillaume Th\'eret}
\address{Guillaume Th\'eret, Lyc\'ee Ni\'epce, 
71100 Chalon-sur-Sa\^one, France\\
\&\ 
Institut de Math\'ematiques de Bourgogne\\
9 avenue Alain Savary,
21078 Dijon, France}
 \email{guillaume.theret71@orange.fr}

\date{\today}

\begin{document}

\begin{abstract}

   The memoir \emph{Theorie der Parallellinien} (1766) by Johann Heinrich Lambert is one of the founding texts of hyperbolic geometry, even though its author's aim was, like many of his predecessors', to prove that such a geometry does not exist. In fact, Lambert developed his theory with the hope of finding a contradiction in a geometry where all the Euclidean axioms are kept except the parallel axiom and that the latter is replaced by its negation. In doing so, he obtained several fundamental results of hyperbolic geometry. This was sixty years before the first writings of Lobachevsky and Bolyai appeared in print.   

\noindent  In the present paper, we present Lambert's main results and we comment on them.
 
 \bigskip 
 
    \noindent        A French translation of the  \emph{Theorie der Parallellinien},  together with an extensive commentary, has just appeared in print \cite{PT}.

\bigskip 


\bigskip 

\bigskip

     \bigskip     
        
\noindent The paper will appear in \emph{Ga\d{n}ita Bh\=ar\=at\=\i } (Indian Mathematics), the Bulletin of the Indian Society for History of Mathematics. 
\bigskip

 \end{abstract}
 
\maketitle

 \section{Introduction} 

The first published treatise on hyperbolic geometry is Lobachevsky's \emph{Elements of geometry}  \cite{Loba-Elements}, printed in installments in the  \emph{Kazan Messenger} in the years 1829-1830. Before that article, Lobachevsky wrote a memoir on the same subject, which he presented on the 12th (Old Style; 23th New Style) of February 1826 to the Physico-mathematical Section of Kazan University. The 
 title of the memoir is {\it Exposition succinte des principes de la g\'eom\'etrie avec une d\'emonstration rigoureuse du th\'eor\`eme des parall\`eles} (A brief Exposition of the principles of geometry with a rigorous proof of the theorem on parallels).  The manuscript of the memoir does not survive; it was ``lost" by the referees. The interested reader may consult the recent edition of Lobachevsky's \emph{Pangeometry} \cite{L} which contains a detailed report on Lobachevsky's published books and articles. 
 
 Before Lobachevsky's works, one can find a few texts which contain results on hyperbolic geometry which their authors developed in an attempt to show that such a geometry does not exist. As a matter of fact, these authors were hoping that in examining a geometry based on the negation of Euclid's parallel axiom,\footnote{Euclid's parallel axiom is sometimes referred to as Euclid's postulate. There is a subtle difference between an ``axiom" and a ``postulate" which can be traced back to the Greek philosophy of science (it is discussed in particular by Aristotle), but we shall not be concerned with this question in the present paper. We shall use both terms, axiom and postulate, without distinction.} one could reach conclusions that would contradict the other Euclidean postulates. The problem of whether the parallel axiom is a consequence of the other Euclidean axioms is the famous \emph{parallel problem}, one of the most important mathematical problems in all history. It is important because the volume of works that were dedicated to it, and because of the conclusion of these efforts, namely, the invention of hyperbolic geometry. We shall consider this problem in some detail in \S \ref{s:parallels} of this paper. 
   
 The memoir \emph{Theorie der Parallellinien} (Theory of parallel lines) by  Johann Heinrich Lambert (1727-1777), written probably in 1766\footnote{The editors of \cite{Engel-Staeckel} state that this date was given by Daniel Bernoulli, who edited and published this work after its author's death.} and which is our subject matter in the present paper,  is a masterpiece of mathematical literature, and its author is one of the most outstanding minds of all times.   Besides being a mathematician, he was an astronomer, physicist (he wrote on optics, magnetism, acoustics, etc.), cartographer, philosopher and linguist. His achievements in mathematics include the first proof that $\pi$ is irrational, and he also conjectured that this number is transcendental.  The book \cite{PT} contains a biography of Lambert, based on the older biographies \cite{Eloge}, \cite{GT}, \cite{Remy} and  \cite{Wolf}. For a lively exposition of some aspects of Lambert's life, we refer the reader to the report \cite{Raghu} by R. Raghunathan published in the present journal.

 Lambert wrote his \emph{Theorie der Parallellinien} in an attempt to prove, by contradiction, the parallel postulate. He deduced remarkable consequences from the negation of that postulate. These consequences make his memoir one of the closest (probably \emph{the} closest) text to hyperbolic geometry, among those that preceded the writings of Lobachevsky, Bolyai and Gauss. We recall by the way that hyperbolic geometry was acknowledged by the mathematical community as a sound geometry only around the year 1866, that is, one hundred years after Lambert wrote his memoir.

 To give the reader a feeling of the wealth of ideas developed in Lambert's memoir, let us review some of the statements of hyperbolic geometry that it contains. Under the negation of Euclid's parallel postulate, and if all the other postulates are untouched, the following properties hold. (The statements are planar): 
\begin{enumerate}
\item The angle sum in an arbitrary triangle is less than $180 \degree$.
\item The area of triangles is proportional to \emph{angle defect}, that is, the difference between $180 \degree$ and the angle sum.
\item There exist two coplanar disjoint lines having a common perpendicular and which diverge from each other on both sides of the perpendicular.
\item Given two lines coplanar $d_1$ and $d_2$ having a common perpendicular, if we elevate in the same plane a perpendicular $d_3$ to $d_1$ at a point which is far enough from the foot of the common perpendicular, then $d_3$ does not meet $d_2$.
\item Suppose we start from a given point in a plane the construction of a regular polygon, putting side by side segments having the same length and making at the junctions equal angles having a certain value between $0$ and $180\degree$ (see Figure \ref{2-Lambert3} below). Then, the set of vertices of these polygons is not necessarily on a circle. Equivalently, the perpendicular bisectors of the segment do not necessarily intersect.
\item \label{5} There exist canonical measures for length and area. 
\end{enumerate}
There are other properties which are discussed in \S \ref{s:commentary} of the present paper.

Property (\ref{5}) may need some comments. There are several ways of seeing the existence of such a canonical measure. For instance, we know that in hyperbolic geometry, there exists a unique equilateral triangle which has a given angle which we can choose in advance (provided it is between $0$ and $60\degree$). This establishes a bijection between the set of angles between $0$ and $60\degree$ and the set of lengths. We know that there is a canonical measure for angles (we take the total angle at each point to be equal to four right angles.) From the above bijection, we deduce a canonical measure for length. This fact is discussed by Lambert in \S 80 of his memoir. Several years after Lambert, Gauss noticed the same fact. In a letter to his friend Gerling, dated April 11, 1816 (cf. C. F. Gauss,  {\it Werke}, Vol. VIII, \cite{Gauss} p. 168), he writes: ``It would have been desirable that Euclidean geometry be not true, because we would have an a priori universal measure. We could use the side of an equilateral triangle with angles  59${}^{\mathrm{o}}$59'59,9999'' as a unit of length". We note by the way that there is also a canonical measure of lengths in spherical geometry, and in fact, a natural distance in this geometry is the so-called ``angular distance".

It also follows from Lambert's memoir that in some precise sense there are exactly three geometries, and that these geometries correspond to the fact that in some (equivalently, in any) triangle the angle sum is respectively equal, greater than, or less than two right angles. This observation by Lambert is at the basis of the analysis that he made of the quadrilaterals that are known as \emph{Lambert quadrilaterals}, or \emph{Ibn al-Haytham--Lambert} quadrilaterals. These are the trirectangular quadrilaterals (that is, quadrilaterals having three right angles), and Lambert studied them systematically, considering successively the cases where the fourth angle is obtuse, right or acute.  It is fair to note here that Lambert was not the first to make such an analysis in the investigation of the parallel problem, and we mention the works of Gerolamo Saccheri (1667-1773) \cite{Saccheri} and, before him, Ab\=u `{}Al\=\i \  Ibn al-Haytham (see Footnote \ref{f:Ibn}) and `Umar al-Khayy\=am (see Footnote \ref{f:Khayyam}). The three geometries suggested by Lambert's and his predecessors' analysis correspond to constant zero, positive or negative curvature respectively, but of course Lambert and his predecessors did not have this notion of curvature. The interpretation of the three geometries in terms of curvature  was given one century after Lambert's work, by Beltrami.

Another important general property on which Lambert made several comments and which he used thoroughly in his memoir is the following: there exist strong analogies between statements in the three geometries (Euclidean, spherical and hyperbolic), with the consequence that some of the statements in the three geometries may be treated in a unified manner. More precisely, he noticed that there exist propositions that are formally identical in the three geometries up to inverting some inequalities or making them equalities. A well-known example is the fact that in Euclidean (respectively spherical, hyperbolic) geometry, the angle sum of triangles is equal to (respectively greater than, smaller than) two right angles. Another example, contained in Lambert's memoir, is that in an equilateral triangle $ABC$, if $F$ is the midpoint of $BC$ and $D$ the intersection point of the medians, we have $DF=\frac{1}{3} AF$ (respectively $DF>\frac{1}{3} AF$, $DF<\frac{1}{3} AF$) in Euclidean (respectively spherical, hyperbolic) geometry. There exist several statements of the same type, in which one passes from one geometry to the other by inverting certain inequalities. Euclidean geometry appears in this setting as the frontier geometry between spherical and Euclidean geometries. In relation with this, Lambert noticed that certain formulae of hyperbolic geometry can be obtained by replacing, in certain formulae of spherical geometry, distances by the same distances multiplied by the imaginary number $\sqrt{-1}$, and by keeping angles untouched. One well-known example is the following: Take, as a model of spherical geometry, the sphere of radius $r$ (or curvature $1/r^2$). Recall that the area of a spherical triangle is equal to $r^2(\alpha+\beta+\gamma-\pi)$, where $\alpha,\beta,\gamma$ are the angles  (in radians). This result is attributed to Albert Girard (1595-1632), who stated it in his  \emph{Invention nouvelle en alg\`ebre}  (1629). If instead of $r$ we take an \emph{imaginary radius} $\sqrt{-1}r$, we obtain, as a formula for area, $-r^2(\alpha+\beta+\gamma-\pi)=r^2(\pi-\alpha-\beta-\gamma)$, which is precisely the area of a triangle of angles $\alpha, \beta,\gamma$  in the hyperbolic space of constant curvature $-1/r^2$. Lambert declares at the occasion of a closely related idea that ``we should almost conclude that the third hypothesis\footnote{The ``third hypothesis", in the memoir of Lambert, is the one where the angle sum in a triangle is less than two right angles. This is the case of hyperbolic geometry.} occurs on an imaginary radius". 

In reviewing some theorems of hyperbolic geometry, Poincar\'e in \cite{Poincare} writes the following: ``The theorems are of course very different from those to which we are accustomed and at first sight they are confusing". Lambert was surely puzzled by the properties in the above list, but he dis not conclude that such a geometry cannot exist.

To end this introduction, let us recall that independently of his memoir on parallels, Lambert was an extremely gifted mathematician. We can quote Hadamard, from his review of the book by St\"ackel and Engel containing Lambert's \emph{Theorie der Parallelinien}, cf. \cite{Hadamard}:
\begin{quote}\small
One cannot but admire, with the authors, the genius, which in some sense is prophetic, of this man who was the first to prove that $\pi$ is irrational, who announced the transcendence of this number and who, at a time where complex numbers were completely foreign to geometers, dared the assertion of the identity between non-Euclidean geometry and the geometry on a sphere of imaginary radius.
\end{quote}

The rest of this paper has two sections. In the first one, we present briefly the problem of parallels. In the second one, we review in some detail the memoir of Lambert, trying to highlight the most important ideas it contains. A French translation of this memoir, together with an extensive commentary, has just appeared \cite{PT}.   An English translation of the first ten pages is contained in the book \cite{Ewald}. V. De Risi  informed us that he is preparing a new edition, with an English translation by A. Laywine, to be published by Birkh\"auser.

\section{On the theory of parallel lines}\label{s:parallels}
The expression ``theory of parallel lines", which is  the title  of Lambert's memoir, usually denotes the attempts made, during a period  which lasted approximately two thousand years, to prove the fifth postulate of Euclid's \emph{Elements}. The history of this major question of our mathematical heritage is very complex, and in fact, it has many facets. Contrary to a widespread idea, it is not true that all the geometers who spent some effort investigating that postulate tried to deduce it from the other Euclidean axioms. Some of them just considered that the postulate is complicated and non-intuitive (in comparison with the other postulates) and therefore, they suggested replacing it by another postulate which would be simpler and more intuitive, such as the existence of translations, the existence of equidistant lines, the existence of homothetical triangles, etc. Considered from this point of view, some of these mathematicians succeeded in their goal. Another group of authors kept this postulate unchanged but they changed the definition of parallel lines so that it becomes simpler. A third group of authors tried to deduce the postulate from the other Euclidean axioms. We know that these attempts were doomed to failure.   Poincar\'e, in an essay on that question in 1891 and which we already mentioned (cf. \cite{Poincare}), wrote: ``One cannot imagine the efforts spent in this unrealistic hope".

 The list of mathematicians who worked on the problem of parallels is impressive. Aristotle, who flourished several decades before Euclid, already discussed this problem in his \emph{Prior Analytics} and \emph{Posterior Analytics}. He declared that there is a vicious circle in the theory of parallels, that is, there is something impossible to prove, ``because of a difficulty in the definition of parallels". We refer the reader to \cite{Heath-Aristotle} for a discussion of Aristotle's writings on the question. Among the known geometers in Europe who worked on the parallel problem, we mention Wallis, d'Alembert, Euler, Lagrange, Clairaut, Legendre and Fourier, and there are many others. It is also well known that the three founders of hyperbolic geometry, Lobachevsky, Bolyai and Gauss, before developing that theory, spent a few years in trying to deduce the parallel axiom from the other axioms.

Let us present, as examples, two short texts on the question of parallels: one from Greek antiquity, and the other from the Arabic period.

The first text is by Proclus.\footnote{\label{f:Proclus} Proclus (c.  412 CE - c. 425 CE) was a neo-Platonic philosopher and mathematician who studied in Alexandria. His written production is impressive, and part of it survives, including his \emph{Commentaries on Book I of Euclid's Elements} \cite{Proclus} and several essays on the mathematics contained in Plato's writings. His work on the parallel problem is mentioned in \S 9 of Lambert's memoir.  In his essay \cite{Proclus}, Proclus gives a proof of the parallel axiom based on the assumption that if two coplanar lines are disjoint, their distance is bounded. It is well known that such an assumption is equivalent to the parallel axiom.} We mention this author because of his importance in the history of geometry, but also  because Lambert, at several places in his memoir, refers to him. The text is extracted from his \emph{Commentaries on Book I of Euclid's Elements} \cite{Proclus}, p. 168ff.):
\begin{quote}\small
This [The fifth postulate] must absolutely be deleted from the postulates; it is a theorem which offers several difficulties, which Ptolemy tried to elucidate in a certain book, and whose proof requires too many definitions and theorems. On the other hand, Euclid presents us the converse of this postulate as a theorem. Some, who are surely deceiving themselves, assessed that they can include it among the postulates, because the certainty of the simultaneous inclination and the meeting of the straight lines is given immediately by the lessening with respect to two right angles. 
But Geminus\footnote{Geminus of Rhodes (c.  110 BCE - 40 BCE).} gave them the right answer when he said that those who are the chiefs of this science warned us that in geometric reasoning we should not at all have any regard to things which are only plausible [...] The fact that if straight lines hereby bend on each other and if the angles decrease is true and unavoidable, whereas the fact that they eventually meet each other when they bend more and more on each other when they are produced enough is plausible but not inescapable, unless some reasoning shows that this fact is true for straight lines. Indeed, the fact that certain lines exist which are asymptotic but non-intersecting seems improbable, but it is nevertheless true, and it has been discovered for other forms of lines. Thus, why is it that what is possible for the latter is not possible for straight lines? \end{quote}

We deduce from this text that Proclus considered the possibility that two straight lines which make, with a third line, and on the same side, two angles whose sum is less than two right angles, may not intersect but are asymptotic. Lambert recalls this fact at the beginning of his memoir. Because of that, Proclus is also considered sometimes as a precursor of hyperbolic geometry.

Now we quote a small text from the Arabic eleventh century mathematician,   `Umar al-Khayy\=am.\footnote{\label{f:Khayyam}  `Umar al-Khayy\=am  (c.  1048 - c. 1131) worked with the notion of \emph{motion}, which had been introduced for the first time by Ibn al-Haytham (cf. Footnote \ref{f:Ibn}) as a primitive notion in a system of axioms of Euclidean geometry, cf. \cite{Rashed1}. In modern terms, a motion is called a \emph{rigid transformation}, or a \emph{congruence}. This notion was used by Hilbert in his axioms of geometry, several centuries after Ibn al-Haytham.}  The reason we mention Khayy\=am is that he studied quadrilaterals which are called \emph{isosceles birectangular}.  These are quadrilaterals which have two right angles adjacent to a common side, with the two sides containing these angles being equal (see Figure \ref{2-Lambert6} below). In his study, Khayy\=am made three hypotheses on these two (equal) angles, namely, they can be right, obtuse, or acute. The three hypotheses are analogous to the three hypotheses that Lambert made for the trirectangular quadrilaterals. We reproduce a short text from Khayy\=am's {\it Commentary on the difficulties of certain postulates in the work of Euclid}, written in 1077 (French translation in Rashed in Vahabzadeh (\cite{Rashed0}, p. 308).  Khayy\=am mentions in this text several works done before him, in particular the work of Ibn al-Haytham\footnote{\label{f:Ibn} Ab\=u  `Al\=\i \ Ibn al-Haytham (c.  965 - c. 1039) is an Arabic mathematician who was known in the Latin world under the name of Alh\=azen. He was also an engineer, physicist, astronomer (a crater of the moon carries his name) and he was called the father of modern optics. His work on optics was known in Europe  thanks to a translation by Gerard of Cremona. We owe to Ibn al-Haytham several treatises on geometry, among which the {\it Book of explanations of the postulates of the book of Euclid on the Elements} and the {\it Book on the resolution of the doubts raised by Book I of Euclid's Elements}. In the first treatise, Ibn al-Haytham defines a \emph{parallel} as the locus of the extremities of a segment which moves perpendicularly to a given line. He argues that this locus is a line. He then tries to prove the parallel postulate by examining trirectangular quadrilaterals, making three hypotheses on the fourth angle: right, acute or obtuse. His analysis is presented in \cite{Youschkevitch} p. 117, Rosenfeld \cite{Rosenfeld} p. 59 and  Pont \cite{Pont} p. 169. With this analysis, Ibn al-Haytham was a direct precursor of Lambert.
  The book by Rashed  \cite{Rashed1} contains a translation of and a commentary on the mathematical works of Ibn al-Haytham.} who, before him, considered the trirectangular quadrilaterals, investigating successively the three hypotheses on the fourth angle. 
 \begin{quote}\small
 [Euclid] introduced an important postulate, and he did not prove it: \emph{If two arbitrary straight lines cutting in two points a straight line make from the same side of it angles [whose sum is] less than a right angle, then they intersect on that side.} [...] I observed that several among those who examined his work (like Heron\footnote{Heron of Alexandria (c. 10 CE - 70 CE) was a Greek mathematician and engineer. In mathematics, he is remembered for his formula that gives the area of a triangle in terms of the side lengths and for his iterative method for finding the square root of two.}  and Eutocius\footnote{Eutocius of Ascalon (c. 480 - c. 540) was a Greek mathematician, mostly known for his commentaries on  Archimedes' treatises and on Apollonius' \emph{Conics}.}) and solved the difficulties it contains did not at all consider the difficulty related to this notion. The moderns, like Al-Kh\=azin, Al-Shann\=\i , Al-Nayr\=\i z\=\i , etc., hardly touched upon its proof. None of them reached a proof which is beyond reproach.  On the contrary, each of them postulated a thing which was not easier to admit. And if there were not so many copies of those books and of persons who owned them and studied them, I would have quoted them here, and I would have explained where there is a postulate and where there is an error, despite the fact that it is not really an easy matter to identity this in these works.
 
 I saw a work by Ab\=u `{}Al\=\i \ Ibn al-Haytham, whose title is  \emph{The resolution of the doubts raised by Book I of Euclid's Elements}, and I had no doubt that he applied himself to that assumption and that he proved it. But when I examined it with delight, I realized that the author wanted this postulate to be at the beginning of the book, together with the other principles, with no need for a proof, that he did disproportionate efforts to attain this goal, that he changed the definition of the parallels, and that he did things which are dazzling, all of them exterior to that art.
 \end{quote}

There is a huge list of papers on the parallel problem. D. M. Y. Sommerville's {\it Bibliography of non-Euclidean geometry} (1911) \cite{Sommerville-Bibliography} contains a chronological catalogue of circa 4000 items, from the beginning of the work on the theory of parallels (IVth century  BCE) until the year 1911. The best modern reference for the history of the parallel problem is the book \cite{Pont} by Pont. We also refer to the exposition in \cite{Gray}.

\section{mathematical commentary on Lambert's \emph{Theorie der Parallellinien}}\label{s:commentary}

 Lambert's memoir is divided into three parts: \S 1 to 11, \S 12 to 26, and \S 27 to 88. The central ideas are the following.
 
 In the first part, the author recalls the problem of parallels, presenting Euclid's eleventh axiom, and the position it occupies among the propositions and the other axioms of the \emph{Elements}. He mentions several difficulties presented by this axiom, quoting commentaries and attempts at proofs by his predecessors. Lambert, who was a fervent reader of classical literature, certainly knew the works of the Greek commentators and their successors on the parallel problem. Furthermore, he was aware of Kl\"ugel's dissertation, written in 1763, which contains a description of  28 attempts to prove the parallel axiom. In particular, Lambert knew about Saccheri's work. It is also good to note that Lambert had probably no intention to publish his manuscript in the state it reached us, which explains the fact that certain historical references (in particular to Saccheri) are missing in that manuscript.

In the second part, Lambert presents some propositions of neutral geometry, that is, the geometry based on the Euclidean axioms from which the parallel axiom has been deleted. One reason for which he works out these propositions is that he thinks that they may be used to prove the parallel axiom.

The third part is the most important part of the memoir. Lambert presents his own approach to prove the parallel axiom. He develops a theory based on the negation of that axiom, hoping that it will lead to a contradiction.

We now elaborate on each of these parts.

At the beginning of his memoir, Lambert recalls that several problems on the question of  parallels arose right at the beginning of geometry, even before Euclid's \emph{Elements}. These problems were a challenge for the major geometers. He notes that the difficulty lies in the eleventh axiom of the \emph{Elements} which, if it were an axiom, ought to be clearer and more self-evident. He says that one has the impression that this axiom is a \emph{proposition} and it needs a proof, an impression which is supported by the fact that the converse of the axiom is a  proposition (Proposition 17 of Book I of the \emph{Elements}) and also by the fact that Euclid, who postponed the use of the parallel axiom until the 29th proposition of Book I, developed a large part of geometry without using this axiom. Similar considerations had already been made by several of Lambert's predecessors,  in particular Proclus, as we recalled above. Lambert also mentions (\S 5) the difficulty, which was also raised many times before him, in the necessity of extending the parallel lines at infinity, which obfuscates the
meaning of teh axiom and makes it less intuitive than the other axioms. He then talks about the 
  \emph{representability} of the parallel axiom. He declares that if this axiom is not a theorem, then one should allow lines which are not \emph{straight}, and which are mutually asymptotic (cf. the remark by Proclus which we already mentioned.)
  
In  \S 4 and 5, Lambert criticizes the approach of Christian Wolff\footnote{Christian Wolff (1679-1754) was a German philosopher and mathematician, disciple of Leibniz. He was an ardent advocate of the logico-deductive method in sciences.} to the theory of parallels presented in the latter's \emph{Anfangsgr\"unde aller Mathematischen Wissenschaften} (1710). He recalls that Wolff defined parallel lines as lines which are co-planar and equidistant, a point of view which had already been adopted by several geometers before him and which did not lead to any breakthrough in the question.

One should recall here that Wolff was a world celebrity at the epoch of Lambert, and was considered as the successor of Leibniz. Therefore, his ideas had a great impact, and for that reason his definition of parallel lines was taken up by several authors after him. Lambert's critic  especially concerned the lack of insight of Wolff who considered that one could resolve the problem by changing the definition of parallelism.

At the beginning of \S 11, Lambert brings up two possibilities in order to solve the difficulties raised by the parallel axiom: either to deduce it, as a theorem, from the other axioms, or to replace it by one or several axioms that are simpler and clearer. He writes, at the end of this introductive part: ``I have no doubt that Euclid had also thought of including his eleventh axiom among the theorems". In the rest of the memoir, Lambert will exploit the first path, and he announces a theory that would solve all difficulties.

 The second part of the memoir starts at  \S 12. In this part, the author's main aim is to present a series of propositions which are valid in neutral geometry.  He starts by recalling that Euclid did not make use of the parallel axiom  until Proposition 29, and that after that proposition, he made extensive use of it. He also notes that several propositions after the 29th are equivalent to the parallel axiom. He mentions, as an example, Proposition 32 which says that the  angle sum in a triangle is equal to two right angles. He also states a proposition saying that a line intersects another line if and only if all the parallels to the first line intersect the second.\footnote{This is not a proposition of the \emph{Elements}.} He then announces that he will prove new propositions that do not use the parallel axiom. 
 
 In fact, Lambert states several propositions and he proves some of them, but for the others, he does not succeed, because contrary to what he thinks, these propositions are not valid in neutral geometry. The ``proofs" that he  gives of these propositions are only attempts at proofs, and Lambert acknowledges after these attempts that his proofs are not conclusive. We shall present some these propositions.

\begin{figure}[!hbp]
 \psfrag{a}{$a$}
\psfrag{b}{$b$}
\psfrag{c}{$c$}
\psfrag{G}{$G$}
\psfrag{d}{$d$}
\psfrag{e}{$e$}
\psfrag{H}{$H$}
\psfrag{f}{$f$}
\psfrag{A}{$A$}
\psfrag{B}{$B$}
\psfrag{C}{$C$}
\psfrag{G}{$G$}
\psfrag{D}{$D$}
\psfrag{E}{$E$}
\psfrag{H}{$H$}
\psfrag{F}{$F$}\centering
\includegraphics[width=.65\linewidth]{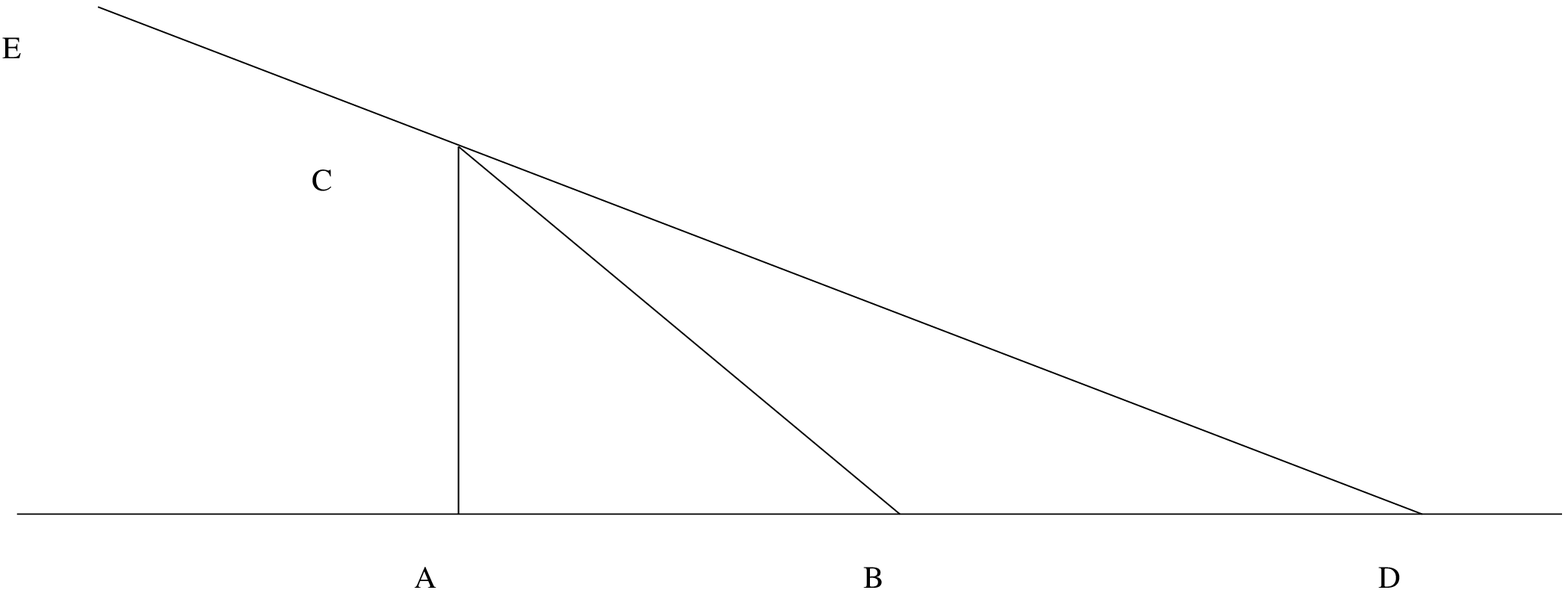}
\caption{\small{}}
\label{2-Lambert1}
\end{figure}

 The first proposition that Lambert proves (\S 13)  says that if  $ABC$ is a triangle with right angle at $A$, and if a line $ED$  passes by $C$ and intersects the line $AB$ at $D$ (Figure \ref{2-Lambert1}), then:
\[\widehat{ACD}< \widehat{ACB}+ \widehat{ABC} < \widehat{ACE}.\]

In \S 15, he proves that if the angle sum in all triangles is constant, then this sum is necessarily equal to two right angles. This is valid in neutral geometry. 
 
Lambert then studies the following construction of regular polygons.

\begin{figure}[!hbp]
 \psfrag{a}{$a$}
\psfrag{b}{$b$}
\psfrag{c}{$c$}
\psfrag{G}{$G$}
\psfrag{d}{$d$}
\psfrag{e}{$e$}
\psfrag{H}{$H$}
\psfrag{f}{$f$}
\psfrag{A}{$A$}
\psfrag{B}{$B$}
\psfrag{C}{$C$}
\psfrag{G}{$G$}
\psfrag{D}{$D$}
\psfrag{E}{$E$}
\psfrag{H}{$H$}
\psfrag{F}{$F$}
\centering
\includegraphics[width=.7\linewidth]{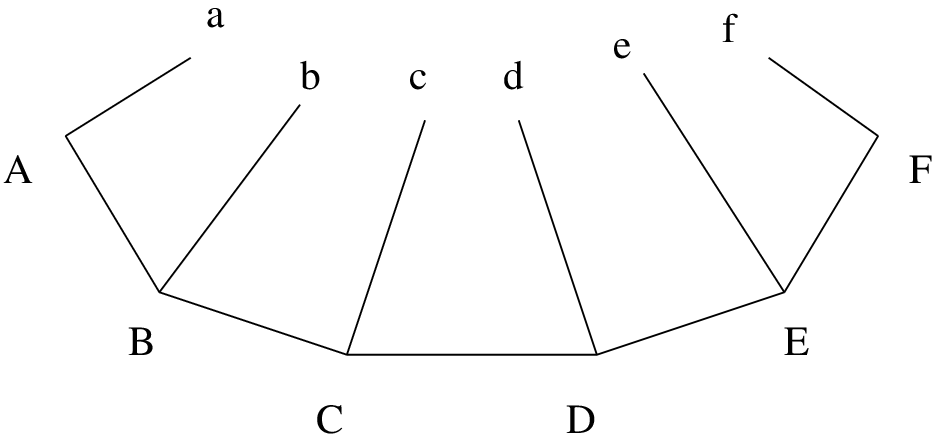}
\caption{\small{}}
\label{2-Lambert3}
\end{figure}

  Take a sequence of points $A,B, C,D,E,F \ldots$ such that  the consecutive segments $AB,  BC, CD, DE, EF, \ldots$ are all congruent\footnote{We use the modern word ``congruent" instead of the Euclidean adjective ``equal".} and the angles at the points  $B,C,D,E, F, \ldots$ are all congruent and strictly less than two right angles
(Figure \ref{2-Lambert3}). 

Consider now the angle bisectors $Aa$, $Bb$, $Cc$, $Dd$,$\ldots$
Lambert states that the points $A,B, C,$  $D,E,$ $F \ldots$ lie on a circle whose center is  the common intersection point of $Aa$, $Bb$, $Cc$, $Dd$, $Ee,\ldots$. It seems that he considered that this proposition is valid in neutral geometry, whereas it is not. It is valid only in Euclidean geometry, and in fact, it is equivalent to the parallel axiom.

\begin{figure}[!hbp]
\psfrag{A}{$A$}
\psfrag{B}{$B$}
\psfrag{C}{$C$}
\psfrag{F}{$F$}
\psfrag{G}{$G$}
\psfrag{D}{$D$}
\psfrag{E}{$E$}
\centering
\includegraphics[width=.70\linewidth]{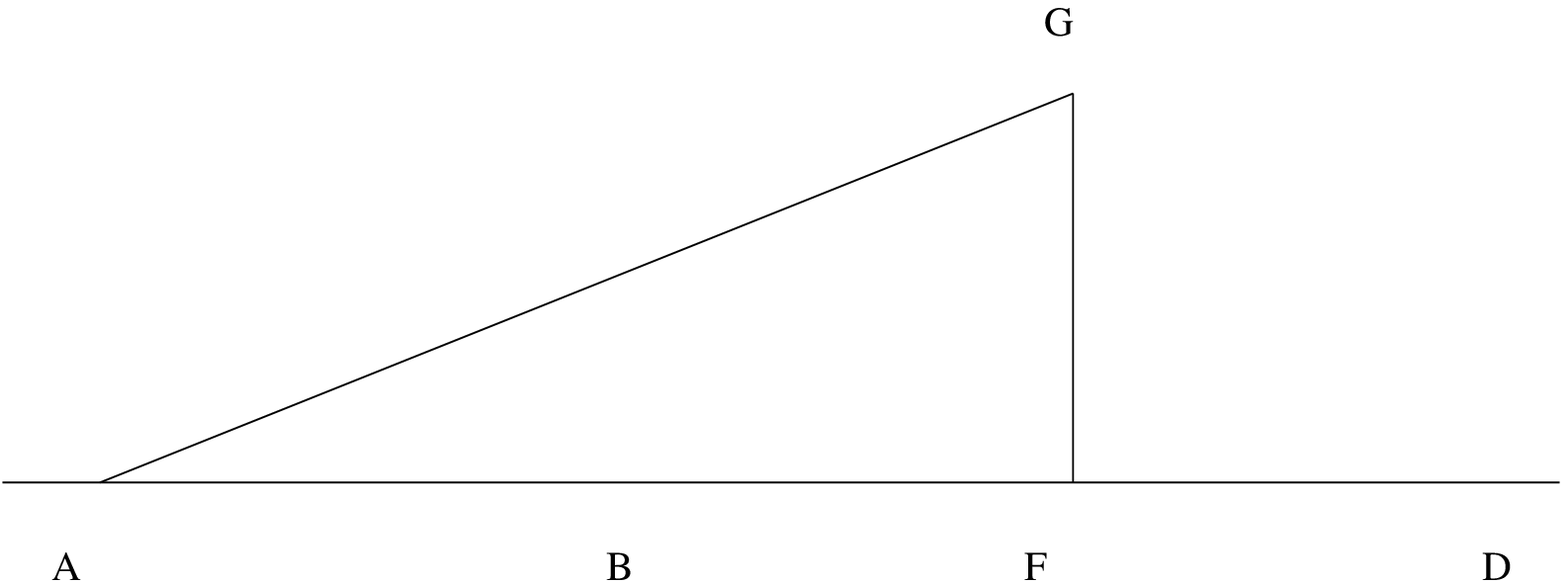}
\caption{\small{}}
\label{2-Lambert-2}
\end{figure}

 In \S 16, Lambert states the following proposition which, like the preceding one, is valid only under the parallels axiom:  Given a line $BD$ and a perpendicular $FG$ elevated at a point $F$ of this line  (Figure \ref{2-Lambert-2}), then we can draw from the point $G$ a line that cuts $BD$ at a point $A$ and makes with it an angle $\widehat{AGF}$ which is as much as we want close to a right angle.

After some preliminaries (\S 17, 18, 19), Lambert sketches three attempts to prove this property in neutral geometry (\S 20, 21, 22), and then acknowledges that he does not succeed.

\begin{figure}[!hbp]
\psfrag{A}{$A$}
\psfrag{B}{$B$}
\psfrag{C}{$C$}
\psfrag{F}{$F$}
\psfrag{p}{$p$}
\psfrag{P}{$\Pi (p)$}
\psfrag{E}{$E$}
\centering
\includegraphics[width=.65\linewidth]{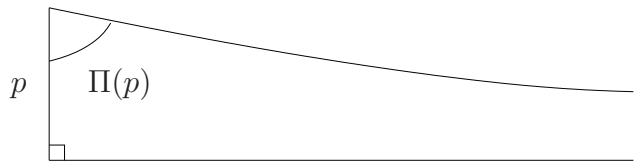}
\caption{\small{The angle $\Pi (p)$ is the angle of parallelism of the segment $p$.}}
\label{loba0}
\end{figure}
Let us take this opportunity to recall that in neutral geometry, the angle $\widehat{AGF}$ has an upper limit which is called  \emph{Lobachevsky's angle of parallelism} (see Figure \ref{loba0}). The value of this angle is a function of the distance from $G$ to $F$ and this angle is right if and only if the parallel axiom is satisfied. In the case of hyperbolic geometry, this upper limit is less than a right angle, cf. \cite{L}.

\begin{figure}[!hbp]
\psfrag{A}{$A$}
\psfrag{B}{$B$}
\psfrag{C}{$C$}
\psfrag{F}{$F$}
\psfrag{G}{$G$}
\psfrag{D}{$D$}
\psfrag{E}{$E$}
\centering
\includegraphics[width=.50\linewidth]{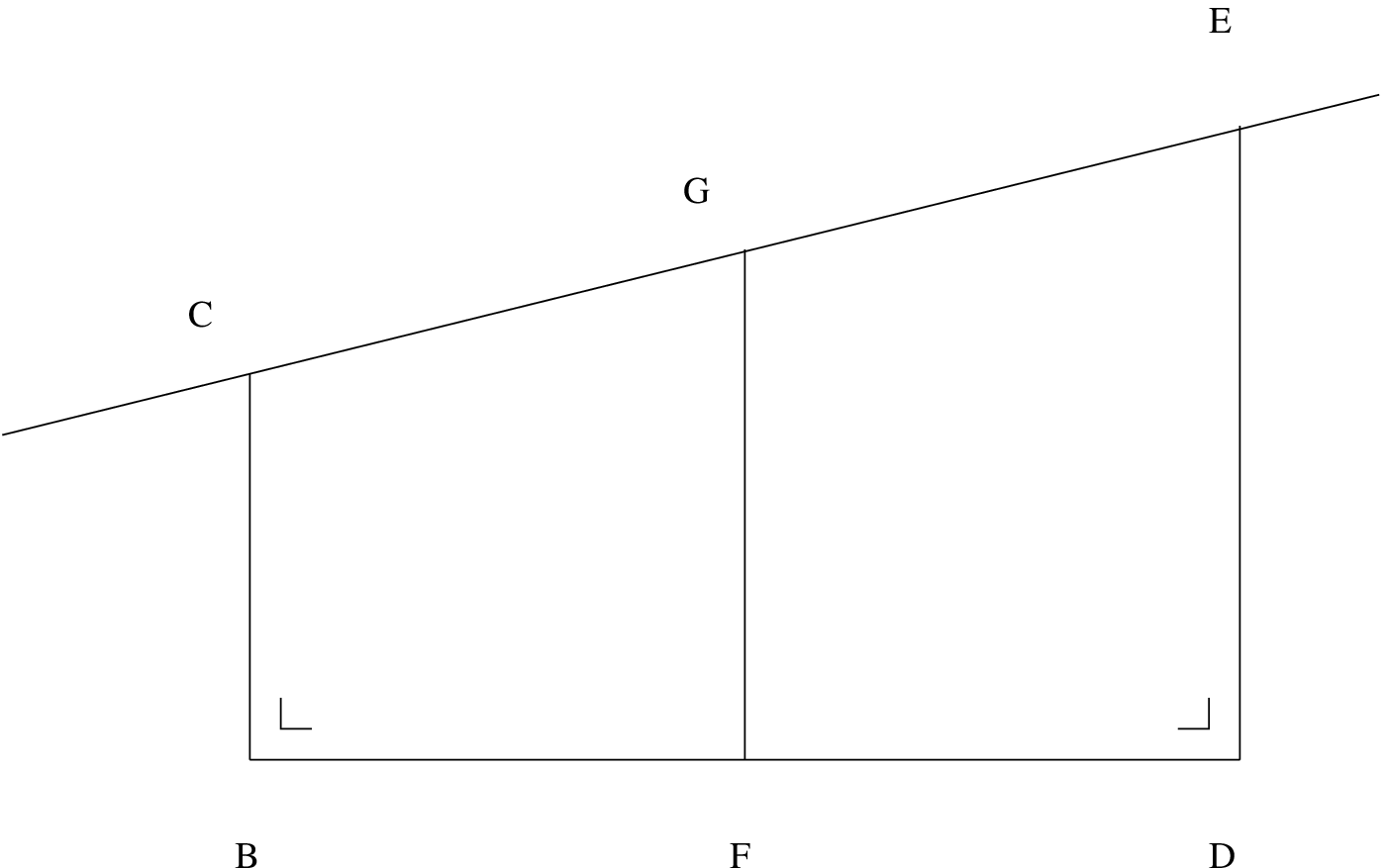}
\caption{\small{}}
\label{3-Lambert2}
\end{figure}

In \S 23, Lambert starts his study of trirectangular quadrilaterals. He considers a quadrilateral $CBDE$ in which the angles $B$ and $D$ are right, and he calls $G$ the intersection point of $EC$ with the perpendicular bisector of $BD$ (Figure \ref{3-Lambert2}). The midpoint of $BD$ is called $F$. He shows that $CB<DE$ if and only if the angle $\widehat{CGF}$ is acute. This is one of several monotonicity properties which are fundamental in the rest of the memoir.
 
\begin{figure}[!hbp]
\psfrag{A}{$A$}
\psfrag{B}{$B$}
\psfrag{C}{$C$}
\psfrag{J}{$l$}
\psfrag{D}{$D$}
\psfrag{L}{$l'$}
\psfrag{H}{$H$}
\psfrag{F}{$F$}
\centering
\includegraphics[width=.60\linewidth]{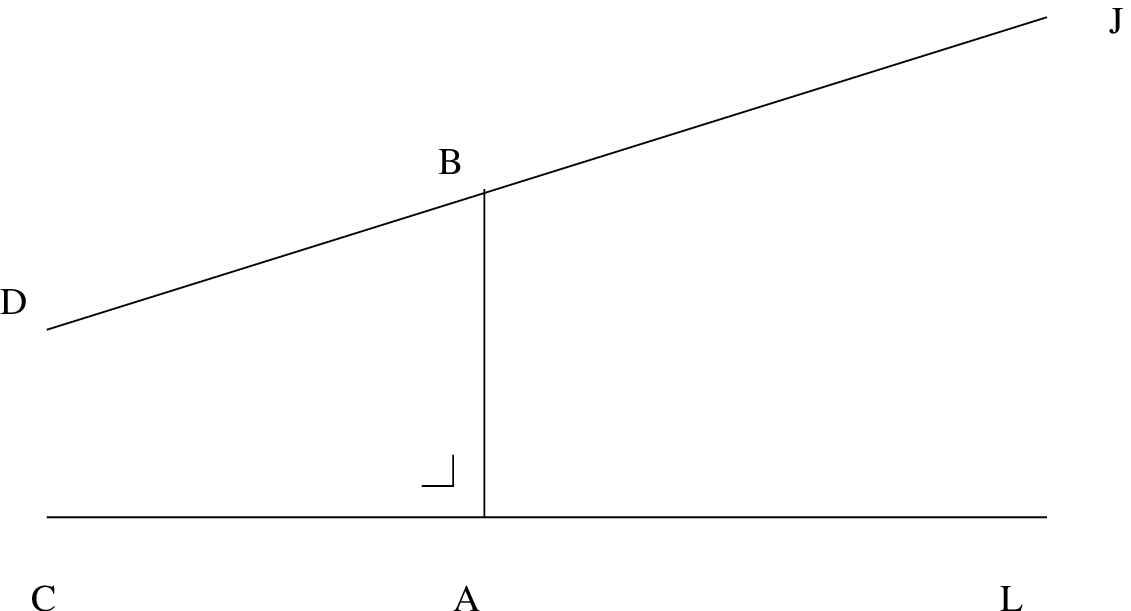}
\caption{\small{}}
\label{2-Lambert5}
\end{figure}

In \S 23, 24 and 25, Lambert studies the distance function  to a line, from a point moving on another line. We know today that this distance function is monotonic only if and only if the parallel postulate holds. Lambert concludes (end of \S 25) that he is not able to prove this property in neutral geometry. 
In \S 26, Lambert proves the following (without using the parallel postulate):
Consider two lines $DJ$ and $CL$ as in Figure \ref{2-Lambert5}, where $B$ is on $DJ$ and $A$ is the foot of the perpendicular drawn from $B$ on the line $CL$. If the angle $\widehat{DBA}$ is always acute (independently of the position of $B$), then this angle is constant.

\bigskip

The third part of the memoir (\S 27 to 88) is the longest. It carries the title ``The theory of parallel lines", the same title as the memoir itself. At the beginning of this part, Lambert recalls that his intention is to present a theory which would fit immediately after Proposition 28 of Book I of Euclid's \emph{Elements}. In other words, this theory would be a sequel to the part of the \emph{Elements} in which the parallel axiom is not used. In this manner, Euclid's Proposition 29 and the subsequent ones, instead of using the parallel axiom, would use Lambert's development. Lambert declares that many statements made by his predecessors, like the one on the equidistance of parallel lines, are not proved, and that his intention is to prove them, after having shown the difficulties which these statements involve. 

He starts by summarizing his approach, where the essential step is to examine successively three hypotheses, namely, that  in a trirectangular quadrilateral, the fourth angle  is either right, or obtuse, or acute. His aim is to show that the second and the third hypotheses contradict the axioms of neutral geometry.

Trirectangular quadrilaterals play a prominent role in the theory of parallel lines. We already mentioned that they have been given later the names \emph{Lambert quadrilaterals}, or \emph{Ibn al-Haytham-Lambert quadrilaterals}.

\begin{figure}[!hbp]
\psfrag{A}{$A$}
\psfrag{B}{$B$}
\psfrag{C}{$C$}
\psfrag{c}{$c$}
\psfrag{D}{$D$}
\psfrag{d}{$d$}
\centering
\includegraphics[width=.5\linewidth]{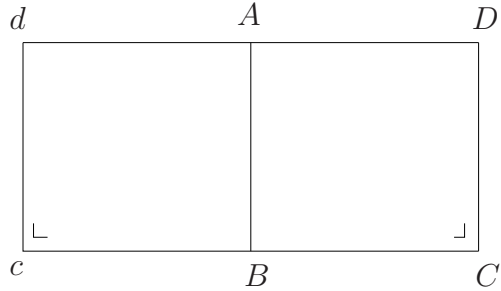}
\caption{\small{The isosceles birectangular quadrilateral $cCDd$ is obtained by folding the trirectangular quadrilateral $ABCD$ along the side $AB$. The angles at $c$ and $C$ are right and $cd=CD$.}}
\label{2-Lambert6}
\end{figure}

In \S 28 to 38, Lambert recalls some properties of parallel lines which are contained in the \emph{Elements}, or which can be derived easily from these properties. He presents some simple properties of trirectangular quadrilaterals. In \S 39, he states again the three hypotheses that concern the  trirectangular quadrilateral $ABCD$ represented in Figure \ref{2-Lambert6}, where the angles $A$, $B$ and $C$ are right. This figure, extracted from Lambert's text, contains at the same time an isosceles birectangular quadrilateral $cCDd$ in which  $cd=CD$, obtained by folding the quadrilateral $ABCD$ along the side $AB$. 

The three hypotheses made by Lambert are then the following:\begin{enumerate}
\item $\widehat{BCD} = 90 \degree$ (\emph{the right angle hypothesis})~;
\item $\widehat{BCD} > 90 \degree$ (\emph{the obtuse angle hypothesis})~;
\item $\widehat{BCD} < 90 \degree$ (\emph{the acute angle hypothesis}).
\end{enumerate}
The first hypothesis is analyzed in \S 40 to 51, the second in  \S 52 to 64, and the third in \S 65 to 88.

Under the first hypothesis, the quadrilateral $cCDd$ is a \emph{rectangle} (that is, a quadrilateral with four right angles, and this implies that the two pairs of opposite sides are congruent). Lambert shows that the first hypothesis implies that for any isosceles birectangular quadrilateral, the perpendiculars drawn from the basis to the opposite side are all congruent (\S 41 and 42). He then proves (\S 43  to 49) that any line $L$ that passes from a point on a side of the quadrilateral making an obtuse angle with that side necessarily intersects the opposite side. The proof  is based on the following monotonicity property: putting side by side along a common basis a sequence of quadrilaterals congruent to $ABCD$, the lengths of the various vertical intersections of the line $L$ with the sides it crosses in the successive quadrilaterals decreases by a quantity which is bounded below. The statement is in  \S 47.  A consequence of this proposition is that the parallel axiom is valid under the hypothesis of the right angle. The aim of the rest of his memoir is to show that the second and third hypotheses are not compatible with the axioms of neutral geometry. 

In \S 48  to 50, Lambert presents the difficulties encountered if we follow the arguments in the reverse direction.

 At the end of this part where he discusses the consequences of the first hypothesis, Lambert states a remarkable result, namely, that the first hypothesis made on some particular trirectangular quadrilateral implies the same property for all trirectangular quadrilaterals.\footnote{Lambert explicitly states such a result only in the case of the first hypothesis. Saccheri, in his article \cite{Saccheri}, had already proved such a ``particular implies general" result (Propositions V, VI and VII of \cite{Saccheri}). Bonola \cite{Bonola} calls \emph{Saccheri's Theorem} the result saying that if in a particular triangle the angle sum is less than (respectively equal, greater than) two right angles, then the same property is valid in any triangle.} In particular, the first hypothesis excludes the second and the third.  

\begin{figure}[!hbp]
\psfrag{A}{$A$}
\psfrag{B}{$B$}
\psfrag{C}{$C$}
\psfrag{D}{$D$}
\psfrag{E}{$E$}
\psfrag{F}{$F$}
\psfrag{G}{$G$}
\psfrag{H}{$H$}
\centering
\includegraphics[width=.7\linewidth]{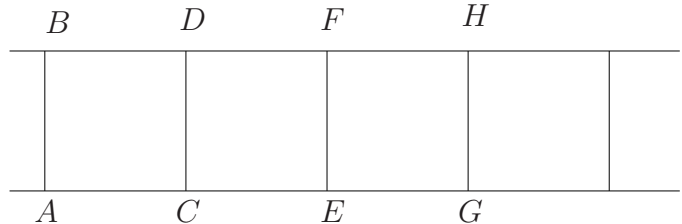}
\caption{\small{At the points $A$ and $B$ and at the points  $C,E,G,...$ the angles are right. Under the second hypothesis, the segments $CD, EF, GH,...$ decrease more than linearly with respect to the distances from the points  $C,E,G,...$ to the point $A$. Likewise, the more the angle at $D,F,H,...$ is far, the more obtuse it is.}}
\label{long-decr}
\end{figure}

In \S 52  to 54, Lambert proves some elementary properties of trirectangular quadrilaterals under the hypothesis of the obtuse angle. In particular, he obtains the fact that any one of the two sides that are bounded by two right angles is longer than the side opposite to it. He then considers (\S 55  to 57) perpendiculars drawn on the line containing a side which is bounded by two right angles, from points situated on the line containing the opposite side and which are outside the rectangle, on the side of the obtuse angle, and whose distances to that obtuse angle are increasing (Figure \ref{long-decr}). 
He shows that the sequence of lengths of these perpendicular segments is monotonically decreasing, and that the difference between a segment and the next one (which is a positive quantity) increases as the sequence of points tends to infinity. Such monotonicity properties are characteristic of non-Euclidean geometry. Analogous properties in spherical geometry are contained in the \emph{Spherics} of Theodosius \cite{Theodosius} (Propositions 6  to 10 of Book III) and  those of Menelaus  (\cite{RP}, Propositions 46 to 61 and 81 to 86).\footnote{Theodosius of Tripoli (c. 160 BCE - c. 100 BCE) and Menelaus of Alexandria (c.  70 CE - c. 140 CE) both wrote treatises on spherical geometry, which they called the \emph{Spherics}. These are the two most important Greek works on spherical geometry that survive, and they are comparable to Euclid's \emph{Elements}.  There exists a French translation  of Theodosius' \emph{Spherics} \cite{Theodosius}.  No Greek manuscript of Menelaus' \emph{Spherics} survives, and  this work reached us only through Arabic translations. There exists a German translation \cite{Krause} based on the Arabic text of Ibn `Ir\=aq. A critical edition with an English translation of that work, together with mathematical commentaries, will appear, based on the Arabic text of  Al-Haraw\=\i \ \cite{RP}. Let us note that the works of Theodosius and Menelaus are very different in spirit and in the methods they use. The methods of Theodosius are based on the Euclidean geometry of the three-dimensional space that is ambient to the sphere, whereas Menelaus works on the sphere itself and his methods are intrinsic. His arguments are mainly based on the fact that this is a geometry in which triangles have angle sum greater than two right angles.}

After this monotonicity for the lengths of segments, Lambert proves (\S 58  to 60) a monotonicity property for angles that the perpendiculars described above make with the line containing the points from which they are drawn. The farthest the angles are, the more obtuse they are (Figure \ref{long-decr} again). Similar monotonicity properties for angles are known to hold  in spherical geometry.\footnote{Propositions 84 to 86 in Menelaus' \emph{Spherics} \cite{RP}.}

Lambert concludes that the lengths of the perpendicular segments cannot be shortened in an asymptotic manner  (\S 62) but that this length becomes zero at some point. Therefore the line which contains the points from which the perpendiculars are produced intersects the line containing the opposite side of the quadrilateral. This contradicts the Euclidean axiom saying that two lines cannot have more than one intersection (\S 64). Lambert concludes (correctly) that the second hypothesis cannot hold.

Lambert understood that his second hypothesis is valid on the sphere (on which some of the other Euclidean axioms not satisfied), but he will only say this in \S 82, during his analysis of the third hypothesis, when he declares (by comparison) that this last hypothesis is realized on a sphere of imaginary radius.

Starting from \S 65, Lambert examines the third hypothesis, the one of the acute angle. He starts by considering a trirectangular quadrilateral (whose fourth angle is therefore acute). He proves (\S 66 and 67) that each of the sides adjacent to the acute angle is greater than the side which is opposite to it. He then shows (\S 68) that each perpendicular drawn from one point on the extension of one side of the acute angle to the opposite side is smaller than a perpendicular drawn from a point on this same line, but further away than the first one. Like in the case of the hypothesis of the obtuse angle, Lambert studies (\S 69) the monotonicity of the angles formed by such lines. The farthest they are, the more acute they are. Then, as in the case studied before, Lambert shows (\S 70) that the lengths of the perpendiculars not only increase, but they can become bigger than any quantity given in advance. Whereas under the second hypothesis such a monotonicity property is shown, without too much difficulty, to contradict one of the axioms of the \emph{Elements}, this is not true in the case of the third hypothesis. At this point, Lambert discovers a property which Saccheri had already highlighted, namely, that under the third hypothesis, there exist disjoint lines which have a common perpendicular and which diverge from each other on each side of that perpendicular.  Lambert is struck by this property, which is very different from what happens in Euclidean geometry. However, he does not conclude that this is impossible.
Let us note that if Lambert had considered, as several of his predecessors, that parallelism implies equidistance or that equidistance defines parallelism, then his aim (finding a contradiction in the negation of the parallel axiom) would have already been attained. But Lambert looks for another argument.
 
\begin{figure}[!hbp]
\psfrag{A}{$A$}
\psfrag{B}{$B$}
\psfrag{C}{$C$}
\psfrag{G}{$G$}
\psfrag{D}{$D$}
\psfrag{E}{$E$}
\psfrag{H}{$H$}
\psfrag{F}{$F$}
\psfrag{J}{$J$}
\psfrag{K}{$K$}
\centering
\includegraphics[width=.85\linewidth]{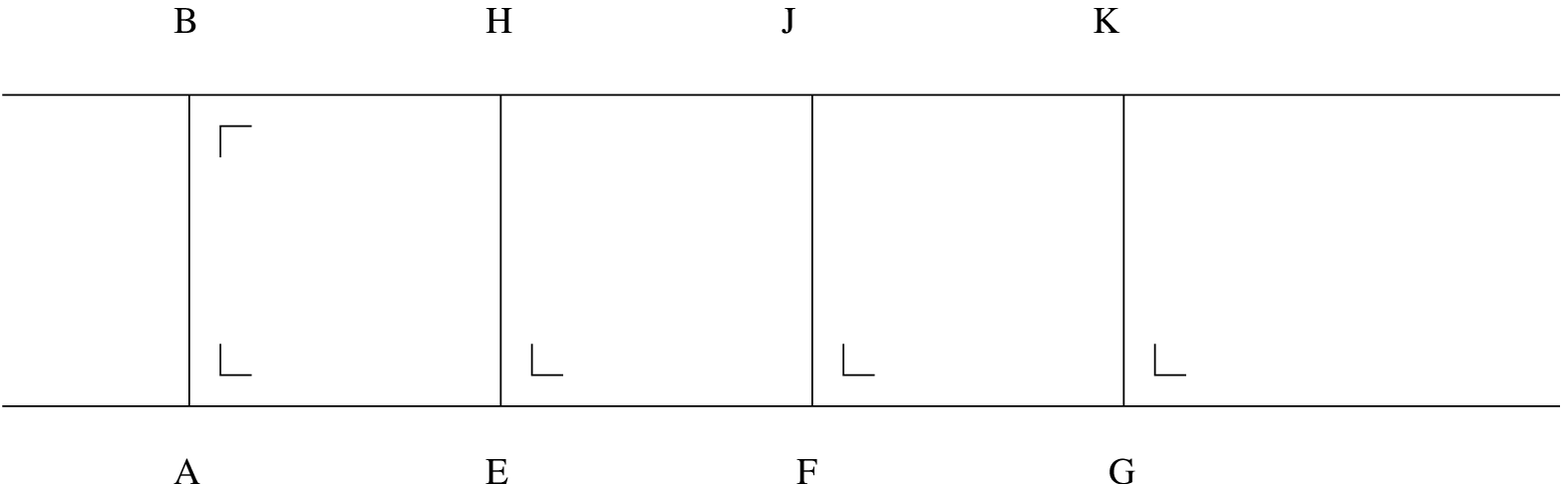}
\caption{\small{}}
\label{2-Lambert7}
\end{figure}

In \S 72, Lambert notices some properties which may occur under the hypothesis of the acute angle. They concern the situation described in Figure \ref{2-Lambert7}, where two lines 
 $AE$ and $BH$ have a common perpendicular $AB$ and where the angles at the sequence of points $E,F,G,\ldots$ are right. The quadrilaterals $BAEH, BAFJ, BAGK,\ldots$ are trirectangular. The first property is that when the points $E,F,G,\ldots$ are far enough from the point $A$, the perpendiculars drawn at these points do not intersect $BH$. The second property is that the angle at $H, J, K,\ldots$ which, as we already saw, decreases when the point moves away from $A$, becomes zero at a certain point, not only as a limiting case, but for a given point on the line $BH$.

In \S 73 and 74, Lambert proves that under the third hypothesis, the angle sum in each triangle is  $< 180 \degree$. 
  
In \S 76 and 77, he shows that in an equilateral triangle $ABC$,  if $F$ is the midpoint of the side $BC$ and if $D$ is the intersection point of the perpendicular bisectors of the sides of    this triangle, then $DF< \frac{1}{3} AF$.\footnote{A more precise version of Lambert's result is contained in the recent paper \cite{PS}.} He notes at the same time that under the hypothesis of the obtuse angle, one has $DF> \frac{1}{3} AF$ (and it is well known that $DF= \frac{1}{3} AF$ under the hypothesis of the right angle).
  
Lambert notes in \S 79 that there are results which are analogous to those that he just mentioned and to several others that he proved before (monotonicity of lengths, of angles, etc.) which are valid under the hypothesis of the acute angle and which become valid under the hypothesis of the obtuse angle if we invert the sense of the monotonicity in the corresponding statements. He notices however that despite this resemblance in the results, the third hypothesis does not lead to a contradiction as easily as the second one.

One of the most important consequences (according to Lambert) of the third hypothesis is the existence of an absolute measure for length, area and volume. To prove this, he reasons as follows  (\S 80). We start from the fact that there is an absolute measure for angles. From here, Lambert notes that there exists an absolute measure for isosceles birectangular quadrilaterals. Such a quadrilateral is completely determined by the value of its acute angle, and the area of such a quadrilateral depends only on the value of that angle. Thus, the value of the angle is an absolute measure for the area of this quadrilateral. Lambert declares that ``this consequence has such a strong attracting force that it easily forces us to wish that the third hypothesis is true!" But he notes that if this were true, then the trigonometric tables ``would be infinitely long". This can be explained as follows. In non-Euclidean geometry, in a right triangle, the ratio of the length of a small side to the hypothenuse is not only a function of the angle that these two sides contain. It also depends on the size of the triangle.  It is in this sense that we can understand the statement that the trigonometry tables of hyperbolic geometry would be infinite (they depend on the size of the triangle). Furthermore, similarity and proportionality of figures totally disappear in such a geometry; figures may only be represented in their real size. 

In \S 81 and 82, Lambert studies the area function of triangles. He shows that we can define the area by \emph{angle defect}, that is, the difference between the angle sum in the triangle and $180 \degree$. He recalls at the same time that under the second hypothesis (the one of the obtuse angle), we can define the area of a triangle as its \emph{angle excess}, that is, the excess of the angle sum to $180 \degree$. Lambert declares that ``one should almost conclude that the third hypothesis occurs on a sphere of imaginary radius".  

This remark on spheres of imaginary radius is interesting for several reasons, and we already commented on that in the introduction.
  One can note here that Lambert, in his \emph{Observations trigonom\'etriques} \cite{Lambert-Ob}, emphasized the analogies between the trigonometric functions sine and cosine and the hyperbolic functions, noting that the hyperbolic cosine can be considered as the cosine of an imaginary arc, and developing a kind of trigonometry in which the sides of a triangle can be considered as imaginary acs.

Several years after Lambert, in the works of Franz Adolph Taurinus (1794-1874), and then in those of  Lobachevsky, similar statements referring to an imaginary sphere were formulated in more explicit forms. Taurinus is the author of a small book which is also called {\it Theorie der Parallellinien}, which he published in 1825 after having received the encouragements of Gauss and under the influence of his uncle, F. K. Schweikart, a law professor at the University of K\"onigsberg who was an amateur mathematician and a friend of Gauss.  The booklet is reproduced in the same volume by St\"ackel and Engel which contains the text of Lambert. In this memoir, Taurinus also attempts a proof of Euclid's parallel axiom. In trying to find a contradiction in the negation of that axiom, Taurinus ended up, like Lambert, developing the bases of non-Euclidean geometry. After the publication fo this work, Taurinus came across Camerer's new edition (in Greek and Latin) of the first 6 books of Euclid's \emph{Elements}, published in Berlin in 1825, which contained a history of the problem of parallels, and he learned there that  similar studies were made before him by Saccheri and Lambert. The next year, he published another work, \emph{Geometriae prima elementa}, in which he developed the first simultaneous expositions of the formulae of the three geometries: spherical, Euclidean and hyperbolic. (Of course, the existence of the latter was for him, like for Lambert, purely hypothetical.) This work contained the fundamental formulae of hyperbolic geometry, obtained by working on a sphere of imaginary radius. Taurinus called such a geometry a {\it  logarithmic-spherical geometry} (\emph{Logarithmisch-sph\"arischen Geometrie}). He noted a passage from the formulae of hyperbolic geometry to those of the spherical which consists in replacing the trigonometric functions  $\sin$ and $\cos$ by the hyperbolic functions $\sinh$ and $\cosh$. In the same manner, Lobachevsky, in his {\it Elements of geometry} (1829-1830) \cite{Loba-Elements}, in his  {\it Geometrical researches on the theory of parallels} (1840) \cite{Loba-Geometrische}, and in other memoirs, highlighted the passage between the formulae of spherical and hyperbolic geometry that consists in replacing, in the formulae of spherical geometry, the side lengths $a,b,c$ of a triangle by the imaginary quantities $a\sqrt{-1},b\sqrt{-1},c\sqrt{-1}$, which of course amounts to replacing the trigonometric functions by the hyperbolic ones. We note by the way that Lambert was among the first mathematicians who systematically used the hyperbolic functions. He is also among the main mathematicians who developed the geometric basis of the theory of complex numbers.

 This passage between the trigonometric formulae of the two non-Euclidean geometries was highlighted by several mathematicians after Lobachevsky, and they were probably attracted by the esthetic appeal of the comparison between the resulting formulae. We mention, as an example,  Beltrami's  {\it Saggio di interpretazione della geometria non-euclidea} (1868) \cite{Beltrami-Saggio}. Beltrami noted there that the trigonometric formulae for the hyperbolic plane  can be obtained from those of the usual sphere by considering the pseudo-sphere (which is a model he had constructed for the hyperbolic plane) as a sphere of imaginary radius $\sqrt{-1}$, and he attributed this observation to E. F. A. Minding  and to D. Codazzi. He writes: 
 \begin{quote}  
 By an observation of Minding (Vol. XX of \emph{Crelle's Journal}), the ordinary formulae for spherical triangles are converted into those for geodesic triangles on a surface of constant negative curvature by inserting the factor $\sqrt{-1}$ in the ratio of the side to radius and leaving the angles unaltered, which amounts to changing the circular functions involving the radius into hyperbolic functions. For example, the first formula of spherical trigonometry
 \[\cos \frac{a}{R} = \cos \frac{b}{R} \cos \frac{c}{R} + \sin \frac{b}{R} \sin \frac{c}{R} \cos A\]
 becomes
 \[\cosh \frac{a}{R} = \cosh \frac{b}{R} \cosh \frac{c}{R} + \sinh \frac{b}{R} \sinh \frac{c}{R} \cos A.\]
\end{quote}
Here, in the case of spherical geometry, $R$ denotes the radius of the sphere, which means that we take the sphere of constant curvature $1/R^{2}$, and in the case of hyperbolic geometry, we consider the plane of constant curvature $-1/R^{2}$.
We note that it is because of such an analogy between the trigonometric formulae for spherical and for hyperbolic geometry that Beltrami chose, for hyperbolic geometry,  the name ``pseudo-spherical geometry". 
     
     We also mention that Klein made, in his {\it \"Uber die sogenannte Nicht-Euklidische Geometrie} (On the so-called hyperbolic geometry) (1871) \cite{Klein-Ueber}, a similar observation: 
 \begin{quote}  
  The trigonometric formulae that hold for our measure result from the formulae of spherical trigonometry by replacing sides by sides divided by $\frac{c}{i}$.
 \end{quote}
     See \cite{ACPK1} for a commentary on Klein's paper.
 
  Closer to us, Coxeter discovered a passage between formulae for volumes of hyperbolic polyhedra and formulae for volumes of spherical ones. In this work, Coxeter made a relation between computations of Lobachevsky in hyperbolic geometry and formulae discovered by Schl\"afli on spherical geometry.

Finally,  Thurston, in his Princeton notes on hyperbolic geometry \cite{Th1},  continues this tradition of using the term ``sphere of imaginary radius" to describe hyperbolic space. He refers to the so-called hyperboloid model of $n$-dimensional hyperbolic geometry embedded in $(n+1)$-dimensional Euclidean space, a model that bears a close analogy with the $n$-sphere embedded in the same Euclidean space. Indeed, whereas the sphere in this model is the space of unit norm vectors for the quadratic form $q(x)= x_0^2+x_1^2+\ldots +x_n^2$, the hyperbolic space is a connected component of the space of vectors of norm $i$ for the quadratic form $q(x)= x_0^2+x_1^2+\ldots -x_n^2$.  In both models, lines (and more-generally, $k$-dimensional planes, for $1\leq k\leq n-1$) are the intersections with the model with planes (or $(k+1)$-planes) passing through the origin, etc.

 Let us conclude our discussion of the \emph{Theorie der Parallellinien}. 

In \S 83, Lambert proves a property which was already mentioned in \S 72. It says (using the notation of Figure \ref{2-Lambert7}) that the perpendiculars drawn at the points
 $E,F,G,\ldots$ on the line $AE$ do not meet the line $BH$, if the point  from which the perpendicular is produced is far enough from $A$.

The last sections of the memoir (\S 84  to 88) contain sketches of some ideas of proofs that the hypothesis of the acute angle may lead to a contradiction, but Lambert acknowledges that all these attempts are fruitless. In the last section, he starts a new attempt. The text is unfinished, and it is not possible to know for sure whether Lambert was convinced whether the parallel postulate is a theorem or not.

\bigskip

{\bf Acknowledgements} 

This work is partially supported by the French ANR project FINSLER. The authors are grateful to the referee and to S. G. Dani for useful comments on the manuscript. The first author wishes to thank Norbert A'Campo for several discussions on non-Euclidean geometry.

 \end{document}